\documentclass[reqno,12pt]{amsart}
\usepackage{amsmath,amsthm,amssymb,tabularx}
\usepackage[cmtip,arrow,all,2cell]{xy}
\usepackage[T1]{fontenc}

\theoremstyle{plain}
\newtheorem{tm}{Theorem}
\newtheorem{lm}[tm]{Lemma}
\newtheorem{prop}[tm]{Proposition}
\newtheorem{cor}[tm]{Corollary}

\theoremstyle{definition}
\newtheorem{remark}[tm]{Remark}

\newcommand{\btm}{\begin{tm}}
\newcommand{\etm}{\end{tm}}
\newcommand{\blm}{\begin{lm}}
\newcommand{\elm}{\end{lm}}
\newcommand{\bpr}{\begin{proof}}
\newcommand{\epr}{\end{proof}}
\newcommand{\beq}{\begin{equation}}
\newcommand{\eeq}{\end{equation}}
\newcommand{\bprop}{\begin{prop}}
\newcommand{\eprop}{\end{prop}}
\newcommand{\bcor}{\begin{cor}}
\newcommand{\ecor}{\end{cor}}

\def\N{\mathbb{N}}
\def\R{\mathbb{R}}

\let\a\alpha
\let\cal\mathcal
\def\n{\mathfrak n}

\def \le {\leqslant}
\def \ge {\geqslant}
\let \<\langle
\let \>\rangle
\let\phi\varphi
\let\kappa\varkappa

\def\liminv{\mathop{\underset\longleftarrow\lim}}
\def\cont{\mathfrak c}

\begin{document}
\author{Yulia Kuznetsova}
\let\thefootnote\relax\footnotetext{Current address: University of Franche-Comt\'e, 25030 Besan\c con Cedex, France. Email: \tt yulia.kuznetsova@univ-fcomte.fr}
\title[Continuity of measurable group representations]{On continuity of measurable group representations and homomorphisms}
\address{University of Luxembourg
6, rue Richard Coudenhove-Kalergi
L-1359 Luxembourg
Tel: (+352) 466644 6802}
%\email{}
\maketitle

\begin{abstract}
Let $G$ be a locally compact group, and let $U$ be its unitary representation on a Hilbert space $H$. Endow the space $\mathcal L(H)$ of linear bounded operators on $H$ with weak operator topology. We prove that if $U$ is a measurable map from $G$ to $\mathcal L(H)$ then it is continuous. This result was known before for separable $H$. We prove also that the following statement is consistent with ZFC: every measurable homomorphism from a locally compact group into any topological group is continuous.
\end{abstract}

Let $G$ be a locally compact group.
We consider its unitary representations, that is, homomorphisms $U$ from $G$ into the group $\mathcal U(H)$ of unitary operators on a Hilbert space $H$. One gets a rich representation theory if the representations considered are weakly continuous, i.e. such that for every $x,y \in H$ the coefficient $f(t)= \<U(t)x,y\>$ is a continuous function on $G$. This requirement is equivalent to strong continuity, i.e. continuity of the function $F(t)=\|U(t)x\|$ for every $x\in H$. Representations satisfying any of these conditions will be further called continuous.

In certain cases it happens that every representation is automatically continuous, as, notably, every finite dimensional unitary representation of a connected semisimple Lie group. This theorem was proved for compact groups by Van der Waerden \cite{waerden} and in general case by A.~I.~Shtern \cite{shtern}. But in general it is easy to construct discontinuous representations, so for automatic continuity, one has to assume some sort of measurability at least. A commonly used notion is as follows. Say that a representation $U$ of a locally compact group $G$ on a Hilbert space $H$ is {\it weakly measurable} if every coefficient $f(t) = \<U(t)x,y\>$ is a measurable function on $G$. Every weakly measurable unitary representation must be continuous if it acts on a separable Hilbert space \cite[Theorem~V.7.3]{gaal}.

However, in general this does not imply continuity: if $G$ is non-discrete, then the regular representation of $G$ on the space $\ell_2(G)$ of countably summable sequences on $G$ is weakly measurable but discontinuous. In this paper we prove that separability restriction can be removed if we use a slightly stronger notion of measurability. Let $\cal L(H)$ be the space of bounded linear operators on the Hilbert space $H$, endowed with the weak operator topology (it is generated by the functions $f_{xy}$ for all $x,y\in H$, where $f_{xy}(A) = \<Ax,y\>$, $A\in \mathcal L(H)$). Say that $U$ is {\it weakly operator measurable} if $U^{-1}(V)$ is measurable for every open set $V\subset \cal L(H)$. Now we can formulate the main result of this paper (Theorem~\ref{rep}): every weak operator measurable unitary representation of a locally compact group is continuous.

The proof is based on a generalization of the so called Four Poles Theorem: if $\cal A$ is a point-finite family of null sets with non-null union in a Polish space, then there is a subfamily in $\cal A$ with a nonmeasurable union (this was proved initially by L.~Bukovsky \cite{bukovsky} and then much simpler by J. Brzuchowski, J. Cicho\'n, E. Grzegorek and C. Ryll-Nardzewski \cite{ryll-1979}). In Lemma \ref{four-poles-continuum}, we prove the same result for subsets of any locally compact group, with a restriction that cardinality of $\cal A$ is not more than continuum.

The second part of the paper deals with automatic continuity of more general group homomorphisms. Most actively this question is studied for homomorphisms between Polish groups, see a recent review of C.~Rosendal \cite{rosendal}. The notion of Haar measurability of $f: G\to H$ is here replaced by universal measurability: the inverse image of every open set is measurable with respect to every Radon measure on $G$.
It is known that every universally measurable homomorphism from a locally compact or abelian Polish group into a Polish group, or from a Polish group to a metric group is continuous. There are also generalizations to other subclasses of Polish groups by S.~Solecki and Rosendal. We omit results on other types of measurability (in the sense of Souslin, Christensen etc.)

If $G$ is not supposed to be Polish, the results are fewer. The most general statement is probably the theorem of A.~Kleppner \cite{kleppner}: every measurable homomorphism between two locally compact groups is continuous. It has been generalized to some special classes of groups by J.~Brzd\k ek~\cite{brzdek}. If one makes no assumptions on the image group, it seems inevitable to impose additional set-theoretic axioms instead. The only result in this direction known to the author, belongs to J.~P.~R.~Christensen \cite{chris}: under Luzin's hypothesis, every Baire, in particular, every Borel measurable homomorphism from a Polish group to any topological group is continuous. Our Theorem \ref{coroll_cont_hom} is proved under Martin's axiom (MA): every measurable homomorphism from a locally compact group to any topological group is continuous.

Theorem \ref{coroll_cont_hom} is reduced to the following question. Let $G$ be a locally compact group; call a set $A\subset G$ {\it extra-measurable} if $SA$ is measurable for any $S\subset G$. An obvious example of an extra-measurable set is an open set. Existence of discontinuous measurable homomorphisms implies existence of null extra-measurable sets; but under MA, as Theorem \ref{nonmeasurable_product} shows, the latter do not exist, so every measurable homomorphism is continuous. In the commutative case, the question of automatic continuity is even equivalent to the existence of a certain sequence of null extra-measurable sets (Theorem \ref{sequence_sets}).

\section{Continuity of unitary representations}

{\bf Definitions and notations}. On a locally compact group $G$, we fix a left Haar measure $\mu$ and the corresponding outer measure $\mu^*$. A map $f:G\to Y$, where $Y$ is a topological space, is called measurable if $f^{-1}(Y)$ is Haar measurable for every open set $U\subset Y$. For a set $A$, $|A|$ denotes its cardinality.

There are two approaches to the construction of Haar measures. One, used by E.~Hewitt and K.~Ross \cite{HR}, yields an outer regular measure: for every measurable set $E$, one has $\mu(E)=\inf\{ \mu(U): E\subset U,\; U \text{ open}\}$. Another one, accepted by D.~Fremlin \cite{fremlin}, leads to an inner regular measure: $\mu(E)=\sup\{ \mu(F): F\subset E,\; F \text{ compact}\}$.

In the $\sigma$-finite case, in particular, on a $\sigma$-compact group, both constructions give the same resulting measure, which is both inner and outer regular. If $G$ is not $\sigma$-compact, the approach of \cite{HR} gives rise to the following pathological sets. A set $A\subset G$ is called {\it locally null} \cite[11.26]{HR} if $\mu(A\cap K)=0$ for every compact set $K\subset G$. Equivalently, $A$ does not contain any set of positive finite measure. Of course, if $A$ is null then $A$ is locally null. Every locally null set $A$ is measurable, and either $\mu(A)=0$ or $\mu(A)=\infty$. In the Fremlin's treatment, there are no locally null non-null sets.

The results of this paper, in particular Theorem \ref{rep}, are valid for both definitions of the Haar measure.

It is known \cite[IV.2.16 and V.7.2]{gaal} that every unitary representation of a locally compact group may be decomposed into a direct sum $U=U_1\oplus U_2$, where $U_1$ is continuous and every coefficient of $U_2$ is (locally) almost everywhere zero. We will say that $U_2$ is {\it singular}. If $U$ acts on a separable space then $U_2=0$ \cite[Theorem V.7.3]{gaal}.
%Proposition V.7.2: if the associated representation of $L_1(G)$ has a cyclic vector then $U_2=0$

Let $U$ act on a Hilbert space $H$. Endow the space $\mathcal L(H)$ of bounded linear operators on $H$ with the weak operator topology (generated by the functions $f_{xy}$ for all $x,y\in H$, where $f_{xy}(A) = \<Ax,y\>$, $A\in \mathcal L(H)$). If $U$ is a measurable map from $G$ to $\mathcal L(H)$, we will say that $U$ is {\it weakly operator measurable}.

The following lemma is known \cite[443P]{fremlin}.

\blm\label{quotient}
Let $G$ be a $\sigma$-compact locally compact group, $K$ its compact normal subgroup, and let $\pi:G\to G/K$ be the quotient  map. If $A\subset G$ is such that $A=AK$ then $A$ is measurable (respectively null) in $G$ if and only if $\pi(A)$ is measurable (null) in $G/K$.
\elm
\iffalse
\bpr The ``if'' part is classical \cite[Corollary~3.3.29]{reiter}. To prove the ``only if'' part, assume that $A=AK$ is measurable. If $\mu(A)<\infty$, its indicator function $I_A$ is in $L_1(G)$. Consider the averaging operator $T:L_1(G)\to L_1(G/K)$, $Tf(\dot x) = \int_K f(x\xi)d\xi$ (in general, $Tf$ is defined almost everywhere on $G/K$, see \cite[\S~3.4]{reiter}). For $\dot x\in G/K$, $TI_A(\dot x) =\int_K I_A(x\xi)d\xi=I_{\pi(A)}(\dot x)$. Thus, $I_{\pi(A)}\in L_1(G/K)$, i.e., $\pi(A)$ is measurable.

Suppose now that $\mu(A)=\infty$. Let $G=\cup_n F_n$, where every $F_n$ is compact. We can assume that $F_n=F_nK$ for all $n$. Denote $A_n=A\cap F_n$, then $A=\cup_n A_n$, $A_n=A_nK$ and $\mu(A_n)<\infty$ for all $n$. As proved above, every $\pi(A_n)$ is measurable, so this is also true for $\pi(A)=\cup_n\pi(A_n)$.
\epr
\fi

The two following facts will be used in further proofs several times, so we prefer to state them separately.

\begin{remark}[Pro-Lie and Polish groups]\label{pro-lie}
Recall that a topological group is called {\it pro-Lie} if it is an inverse (projective) limit of (finite-dimensional) Lie groups (see \cite{hofmor}). It is known that in every locally compact group there is an open pro-Lie subgroup \cite[p.~165]{hofmor}.
If $G$ is a locally compact group and $G=\liminv_{i\in I} G_i$, where every $G_i$ is a Lie group, then these groups can be chosen as $G_i=G/K_i$, where every $K_i$ is a compact normal subgroup of $G$, and the order on $I$ is just inclusion of $K_i$. Every $\sigma$-compact Lie group is Polish (it is first countable, hence metrizable \cite[Theorem~A4.16]{hofmor-compact}, and further apply \cite[Chapter IX, \S~6]{bourbaki}, Corollary of Proposition 2). %If $I$ is countable, $G$ is also Polish (\cite[\S~6, Proposition~1a,b]{bourbaki}).
If all $G_i$ are $\sigma$-compact and $I$ is countable, $G$ is Polish too (\cite[\S~6, Proposition~1a,b]{bourbaki}).
\end{remark}

\begin{remark}[Baire sets in direct products]\label{baire}
Baire sets (\cite[\S51]{halmos}) are the elements of the $\sigma$-algebra generated by all compact $G_\delta$-sets. In the case of a $\sigma$-compact locally compact group (the only case we will need), this is also the $\sigma$-algebra generated by all zero sets of continuous functions (in \cite[11.1]{HR}, this latter property is taken as a definition). Consider the direct product of a family of locally compact groups $\cal G=\prod_{j\in J}G_j$. This is a topological group, which is not necessarily locally compact. Let $G\subset \cal G$ be a closed $\sigma$-compact subgroup. For any $I\subset J$ let $\pi_I:\cal G\to\prod_{j\in I}G_j$ be the natural projection. We say that a set $X\subset G$ depends on coordinates $I\subset J$ if $X=G\cap \pi_I^{-1}(\pi_I X)$. If $F=\cap\, U_n$ is a compact $G_\delta$ set in $G$, then for every $n$ the open set $U_n$ can be chosen as a finite union of basic neighbourhoods in $\cal G$, which depend on finite number of coordinates. It follows that every such $F$, and as a consequence every Baire set depends on a countable set of coordinates.
\end{remark}

\blm\label{four-poles-continuum}
Let $\cal A=\{A_s:s\in S\}$ be a point finite family of null subsets of a $\sigma$-compact locally compact group $G$. If $|\cal A|\le\cont$ and $\cup\cal A$ is non-null, then there is $\cal B\subset \cal A$ such that $\cup \cal B$ is nonmeasurable.
\elm
\bpr
Let $H\subset G$ be an open pro-Lie subgroup, which we can assume compactly generated. Since $G$ is a countable union of $H$-cosets, there is $t_0\in G$ such that $(\cup \cal A)\cap t_0H$ is non-null.
Define $A_s'=(t_0^{-1}A_s)\cap H$ for all $s$, then the family $\cal A'=\{A'_s: s\in S\}$ satisfies all conditions of the theorem and is contained in $H$. Moreover, if a union $\cup\{A'_s: s\in T\}$ is nonmeasurable, then so is $\cup\{A_s:s\in T\}$. Thus, we can assume that $G=H$, i.e. $G=\liminv_{i\in I} G_i$ is $\sigma$-compact and pro-Lie. In this case every $G_i=G/K_i$ is a $\sigma$-compact Lie group, hence Polish. We assume that every $G_i$ is nontrivial, otherwise $G=\{1\}$ and the family $\cal A$ would not exist. %in terminology of Bourbaki, is countable at infinity

We can assume that $S\subset \R$. Let $\mathbb Q=\{q_m:m\in \N\}$ be an enumeration of the rationals, and let $W_{mn}=\cup\{A_s: |s-q_m|<1/n\}$. For every $s\in S$, choose sequences $n_k^{\!(s)},m_k^{\!(s)}$ so that $|s-q_{m_k^{\!(s)}}|<1/n_k^{(s)}$ and $n_k^{\!(s)}\to\infty$ while $n_{k+1}^{\!(s)}>n_k^{\!(s)}$ for all $k$. %$|q_{m_k^{(s)}} - s|$
Then $A_s=\cap_k W_{n_k^{\!(s)} m_k^{\!(s)}}$ for every $s$. Indeed, every point $x\in A_s$ is contained in this intersection; if $x\notin A_s$ then $x\in A_{t_i}$ for at most finite set of points $t_i\ne s$; and every $t_i$ can be separated from $s$ by some interval $|q_{m_k^{(s)}} - s|<1/n_k^{(s)}$, so that $x\notin W_{n_k^{\!(s)} m_k^{\!(s)}}$.

If one of the sets $W_{mn}$ is nonmeasurable, the lemma is proved. Suppose that every $W_{mn}$ is measurable. By \cite[19.30b]{HR}, there exists a Baire set $B_{mn}\subset W_{mn}$ such that $N_{mn}=W_{mn}\setminus B_{mn}$ is null. Further, for every $n,m$ there is a null Baire set $N'_{mn}\supset N_{mn}$. Let $N=\cup_{m,n} N'_{mn}$. Then $N$ is a Baire set, so $W_{mn}\setminus N=B_{mn}\setminus N$ is Baire for all ${m,n}$. Let $W_{mn}\setminus N$ depend on the countable set of coordinates $J_{mn}$. Then every $A_s\setminus N=\mathop{\cap}\limits_k \,( W_{n_k^{\!(s)}m_k^{\!(s)}}\setminus N)$ depends on coordinates $J=\cup_{mn}J_{mn}$, and the set $J$ is countable.

Extending $J$, if necessary, we can assume that the family $\{K_j:j\in J\}$ is closed under finite intersections.  Denote $K=\cap_{j\in J}K_j$. Then $G/K=\liminv_{j\in J}G/K_j$ is a Polish group. Let $\pi:G\to G/K$ be the quotient map, and put $A'_s=\pi(A_s\setminus N)$. Then, since $A_s\setminus N=(A_s\setminus N)K$, the family $\cal A'=\{A'_s: s\in S\}$ is point finite, and by Lemma \ref{quotient} we have that $A'_s$ is null for all $s$, while $\cup\cal A'= \pi(\cup\cal A)$ is non-null. By the Four Poles Theorem for the Polish case \cite{ryll-1979} we get $\cal B'\subset \cal A'$ such that $\cup \cal B'$ is nonmeasurable. Put $\mathcal B=\{A_s:A'_s\in\mathcal B'\}$. Then $\cup\mathcal B\setminus N=\pi^{-1}(\cup\cal B')$ is nonmeasurable, so $\cal B$ is as desired.
\epr

A simple example shows that in the Hewitt\&Ross approach, this theorem is not true for a non-$\sigma$-compact group. Let $\R_d$ be the real line with the discrete topology, and consider the direct product $\R_d\times\R$. Then $X=\R_d\times\{0\}$ is measurable of infinite measure (this is an example of a locally null, non-null set \cite[11.33]{HR}). Every uncountable subset of $X$ is also measurable with infinite measure, and every countable subset is null. Thus, if we put $A_t=\{(t,0)\}$ and $\cal A=\{ A_t: t\in \R_d\}$, then every $A_t$ is null, $\cup\cal A$ is non-null but every subfamily of $\cal A$ has a measurable union. In the Fremlin's approach, this example does not appear.

\btm\label{rep}
Let $G$ be a locally compact group. Then every weakly operator measurable unitary representation of $G$ is continuous.
\etm
\bpr
Let $U:G\to\mathcal L(H)$ be a representation acting on a Hilbert space $H$. Clearly $U$ is continuous if and only if its restriction to any open subgroup is continuous. In $G$, there is an open compactly generated pro-Lie subgroup (e.g., the intersection of an open pro-Lie subgroup and the subgroup generated by a pre-compact neighbourhood of identity). So we can assume that $G=\liminv_{i\in I} G/K_i$ is itself compactly generated and pro-Lie; in particular, $G$ is $\sigma$-compact.

Take any $x\in H$, $\|x\|=1$. Put $f(t)=\<U(t)x,x\>$ and $S=\{t\in G: f(t)\ne0\}$. We can assume that $U$ is singular, then $S$ is null. Suppose towards a contradiction that $U\not\equiv0$, then $e\in S$.

Let us show that $S$ has a null projection onto a Polish quotient of $G$. By \cite[19.30b]{HR}, there exists a null Baire set $B\supset S$. Every Baire set (Remark \ref{baire}) depends on a countable number of coordinates. Let $J\subset I$ be a countable set such that $B=\pi_J^{-1}\pi_J B$. By extending $J$ if necessary, we can assume that the family $\{K_j:j\in J\}$ is closed under finite intersections.  Denote $K=\cap_{j\in J}K_j$. Then $G/K=\liminv_{j\in J}G/K_j$ is a Polish group. Let $\pi:G\to G/K$ be the quotient map and let $S'=\pi(S)$, then $S'\subset\pi(B)$ is null.

Choose an enumeration (probably with repetitions) $P_\a:\a<\cont$ of perfect non-null sets in $G/K$. It is known that there is at most continuum such sets. By induction, we will choose points $t_\a:\a<\n$ with some ordinal $\n\le\cont$ so that $\cup\{t_\a S': \a<\n\}$ is non-null (in $G/K$) and $t_\a\notin \cup\{ t_\beta S' : \beta<\a\}$ for every $\a>0$. Set $t_0=e$. For every $\a$, let $T_\a=\{t_\beta:\beta<\a\}$. If $T_\a S'$ is non-null, stop the procedure. Otherwise $P_\a\setminus T_\a S'\ne\emptyset$, and choose $t_\a$ as any point of this set. Let $\n$ be the ordinal on which we stopped the induction, or $\n=\cont$ if it was not stopped. Set $T=T_\n$. If $\n<\cont$ then as assumed $\mu^*(TS')>0$; if $\n=\cont$ then $TS'$ intersects every non-null perfect set in $G/K$, so it is of full measure. In either case $TS'$ is non-null.

For every $\a<\n$, choose any $z_\a\in \pi^{-1}(t_\a)$ and set $Z=\{z_\a:\a<\n\}$. It follows that $ZSK= \pi^{-1}(\pi(ZS))= \pi^{-1}(TS')$ is non-null in $G$. Recall that $K$ is a normal subgroup, so $ZSK=ZKS$. Define now
\beq\label{formula_Sn}
S_n=\{ t\in G: |f(t)|=|\<U(t)x,x\>|>1/n\}.
\eeq
Then $S=\cup_n S_n$ and $ZKS=\cup_n ZKS_n$. It follows that $ZKS_N$ is non-null for some $N\in\N$.

We claim that the family $\cal A=\{z_\a KS_N: \a<\n\}$ is point-finite. To prove this, first show that $U(z_\a k_2)x\perp U(z_\beta k_1)x$ for any $\a\ne\beta$ and any $k_1,k_2\in K$. Suppose that $\a>\beta$. Then $(z_\beta k_1)^{-1}z_\a k_2\notin S$, because otherwise we would have $z_\a\in z_\beta k_1 S k_2^{-1} \subset z_\beta KSK = z_\beta SK^2=z_\beta SK$ and hence $t_\a=\pi(z_\a)\in \pi(z_\beta SK)=t_\beta S'$, what is impossible by the choice of $t_\a$. This gives us
$$
0=f((z_\beta k_1)^{-1}z_\a k_2) = \<U((z_\beta k_1)^{-1}z_\a k_2)x,x\> = \<U(z_\a k_2)x,U(z_\beta k_1)x\>,
$$
that is, $U(z_\a k_2)x\perp U(z_\beta k_1)x$.

Next, if $t\in z_\a KS_N$ then there is $k_\a\in K$ such that $(z_a k_a)^{-1}t\in S_N$, i.e. $|\< U(t) x, U(z_\a k_\a)x\>| >1/N$. As we have shown above, $U(z_\a k_\a)x$ are orthogonal for different $\a$; since $U$ is unitary, they have norm 1. By Bessel's inequality we have for any $t\in G$:
\begin{align*}
1=\|x\|^2=\|U(t)x\|^2 &\ge \sum_{\a:\, t\in z_\a KS_N} |\<U(t)x, U(z_\a k_\a)x\>|^2
\\&> N^{-2}\cdot |\{\a: t\in z_\a KS_N\}|.
\end{align*}

So $\cal A$ is point finite. It has cardinality $|\cal A|=|Z|=\n\le\cont$ and a non-null union $\cup\cal A = ZKS_N$. Every $z_\a KS_N \subset z_\a KS = z_\a \pi^{-1}(S')$ is a null set. Thus, we can apply Lemma \ref{four-poles-continuum} to get $\mathcal B\subset \mathcal A$ such that $\cup\cal B$ is nonmeasurable. Now recall that $S_N$, by formula \eqref{formula_Sn}, is the inverse image of an open set in $\mathcal L(H)$. The same is true for every translate of $S_N$ and for unions of such translates, in particular for every $z_\a KS_N$ and for $\cup \cal B$. As the inverse image of an open set, $\cup \cal B$ must be measurable. This contradiction proves the theorem.
\epr

\section{Continuity of group homomorphisms}

The contents of this section is valid for both treatments of the Haar measure. For the inner regular variant accepted by \cite{fremlin}, it suffices to ignore the bracketed ``locally'' everywhere.

Let $G$ be a locally compact group. Call a set $A\subset G$ {\it extra-measurable} if $SA$ is measurable for every set $S\subset G$. Every open set is extra-measurable, while a one-point set is not, unless the group is discrete. As shown below (Theorem \ref{cont_hom}), existence of discontinuous measurable homomorphisms implies existence of [locally] null (definition below) extra-measurable sets. It is consistent with ZFC (Theorem \ref{nonmeasurable_product}) that a nonempty [locally] null set cannot be extra-measurable, so it is consistent that every measurable homomorphism from a locally compact group to any topological group is continuous. It is an open question, whether this statement is true in ZFC without any additional axioms. Already in the basic case of the real line the answer is unknown, but for commutative groups one can make the question more precise (Proposition \ref{sequence_sets}).

\btm\label{cont_hom}
Let $G$ be a locally compact group. If there exists a homomorphism $\phi: G\to H$ to a topological group $H$ which is measurable but discontinuous, then there is a family $\cal A$ of nonempty [locally] null extra-measurable sets such that for every $A\in\cal A$:
\begin{align}\label{nbhds}
 &A^{-1}=A;\\
 &\exists B\in\cal A:\quad B^2\subset A;\notag\\
 &\forall x\in G\quad \exists C\in\cal A:\quad x^{-1}Cx\subset A.\notag
\end{align}
\etm
\bpr
Suppose that such $\phi$ exists. Let $U$ be an open neighbourhood of identity in $H$ and let $A=\phi^{-1}(U)$. Then for any $S\subset G$ we have $SA=\phi^{-1} (\phi(S)U)$. Indeed, the inclusion $\phi(SA)\subset \phi(S)\phi(A)=\phi(S)U$ is obvious. For the opposite inclusion, take $z\in \phi^{-1}(\phi(S)U)$ and choose $s\in S$, $a\in A$ such that $\phi(z)=\phi(s)\phi(a)=\phi(sa)$. Let $(sa)^{-1}z=t$, then $t\in\ker \phi$; since $\phi(at)=\phi(a)\in U$, we have $at\in A$ and $z=sat\in SA$.

Now $SA$ is the inverse image of an open set $\phi(S)U$, so it must be measurable. Thus, $A$ is extra-measurable.

Suppose that $\phi^{-1}(U)$ is not [locally] null for every $U$. Take an open neighbourhood of identity $V$ such that $V^{-1}V\subset U$. Then $B=\phi^{-1}(V)$ is by assumption also non-[locally] null. It contains then a set $C$ with $0<\mu(C)<\infty$, so $C^{-1}C$ contains a neighbourhood of identity $W\subset G$ \cite[20.17]{HR}. Then $\phi(W)\subset \phi(C^{-1}C)\subset V^{-1}V\subset U$, so $\phi^{-1}(U)\supset W$. Since $U$ was arbitrary, this implies that $\phi$ is continuous.

Thus, in assumptions of the theorem there is $U$ such that $\phi^{-1}(U)$ is [locally] null. Let $\cal V$ be a base of neighbourhoods of identity in $H$, such that $V\subset U$ and $V=V^{-1}$ for every $V\in \cal V$. Denote $\cal A=\{\phi^{-1}(V): V\in\cal V\}$, then this family has properties \eqref{nbhds}.
\epr

The properties \eqref{nbhds} guarantee that if we take $\cal A$ as a base of neighbourhoods of identity in $G$, this turns $G$ into a topological group \cite[IV,\ \S2]{bourbaki}. However, it does not follow immediately that the identical map on $G$ is measurable, i.e. in general we do not get a converse of this theorem. An equivalence holds in the commutative case:

\bprop\label{sequence_sets}
Let $G$ be a commutative locally compact group. The following are equivalent:\\
(i) There is a homomorphism $\phi: G\to H$ to a topological group $H$ which is measurable but discontinuous;\\
(ii) There is a sequence of [locally] null extra-measurable sets $A_n$ such that for every $n$:\\
 $A_n^{-1}=A_n$; $A_{n+1}^2\subset A_n$.
\eprop
\bpr
(i)$\Rightarrow$(ii): Proved in Theorem \ref{cont_hom}.

(ii)$\Rightarrow$(i): Take the sets $A_n$ as a base of neighborhoods of identity in $G$, then this turns $G$ into a topological group which we can denote $H$. (Note that $H$ is metrizable if $\cap A_n=\emptyset$.) The identical map $\phi: G\to H$ is obviously discontinuous. At the same time, for every open set $U\subset H$ we have $U=\cup T_n A_n$ for some sets $T_n$; every $T_nA_n$ and hence $U=\phi^{-1}(U)$ are measurable, so $\phi$ is a measurable map and (i) holds.
\epr

Existence of sets as in Proposition \ref{sequence_sets}(ii) is an open question even on the real line. Known results on automatic continuity mostly concern Polish groups; here they do not give a ready answer, since the group $H$ obtained in the proof may be not complete (i.e. not Polish).

Let ${\rm add}(\cal N)$ be the minimal cardinality of a family $\mathcal J$ of null sets on the real line $G$ such that $\cup \mathcal J$ is non null. This is called the additivity of the ideal $\cal N$ of Lebesgue null sets in $G$. It is known that additivity of the ideal of Haar null sets is the same for every non-discrete locally compact Polish group \cite[522Va]{fremlin}. It is consistent with ZFC that ${\rm add}(\cal N)<\cont$, but it follows from Martin's axiom (MA) that ${\rm add}(\cal N)=\cont$ (see \cite{fremlin-martin}). This is, in fact, the assumption that we use in our proof. It is known that Martin's axiom follows from the Continuum hypothesis, but is consistent also with its negation. For further discussion of Martin's axiom, we refer to the Fremlin's monograph \cite{fremlin-martin}.

A. Kharazishvili has indicated in private correspondence that the following statement holds for a commutative Polish group:

\blm[MA]\label{polish_group}
Let $G$ be a locally compact Polish group, and let $A\subset G$ be a nonempty set of measure zero. Then there is a set $S\subset G$ such that $SA$ is nonmeasurable.
\elm
\bpr
If $G$ is countable, then by local compactness it has isolated points, and the measure of every point is positive. Then the set $A$ in assumption cannot exist. Thus, $G$ is uncountable without isolated points. Note that $G$ is $\sigma$-compact \cite[Theorems 3.3.1, 3.8.1, 3.8.C(b)]{engelking}.

We will construct $S$ so that both $SA$ and $G\setminus SA$ intersect every perfect set of positive measure. Then $SA$  must be nonmeasurable, since the inner measure of $SA$ and $G\setminus SA$ is zero.

By translating $A$, and then $S$, if necessary, one can assume that $e\in A$. Note that $A^{-1}$ has also measure zero---this follows, e.g., from \cite[20.2]{HR} or \cite[442K]{fremlin}.
%$$
%\mu(A^{-1})=\int I_{A^{-1}}(x)dx = \int I_A(x^{-1})dx = \int I_A(x)\Delta(x^{-1})dx = \int_A \Delta(x^{-1})dx = 0,
%$$
%where $\Delta$ is the modular function of $G$.

Since $G$ is separable and uncountable, there is exactly continuum many closed sets in it. Let $\{P_\xi: \xi<\cont\}$ be an enumeration of all perfect sets of positive measure. By induction, we will choose $s_\xi,d_\xi \in P_\xi$ so that the condition $d_\xi \in P_\xi\setminus s_\eta A$ holds for every $\xi,\eta$. Then $S=\{s_\xi: \xi<\cont\}$ will be as needed, since $s_\xi\in S\cap P_\xi\subset SA\cap P_\xi$ and $d_\xi\in P_\xi\!\setminus \!SA$, so both $P_\xi\cap (SA)$ and $P_\xi\cap(G\setminus SA)$ are nonempty.

Suppose that for all $\eta<\xi$ such $s_\eta,d_\eta$ are chosen, or that $\xi=0$ (the base of induction). Set $D_\xi=\{d_\eta:\eta<\xi\}$ and note that $|D_\xi|<\cont$. Since $P_\xi$ cannot be covered by a less that continuum translates of $A^{-1}$ (here we use the Martin's axiom), we can choose a point $s_\xi\in P_\xi\setminus D_\xi A^{-1}\ne\emptyset$. This implies $(s_\xi A)\cap D_\xi=\emptyset$.

Next, set $S_\xi=\{s_\eta:\eta\le\xi\}$. Then $|S_\xi|<\cont$ and similarly we can choose $d_\xi\in P_\xi \setminus S_\xi A$. By this choice, we have $d_\xi\notin s_\eta A$ for all $\eta\le\xi$, and for $\eta>\xi$ we have $d_\xi\notin s_\eta A$ by the choice of $s_\eta$. This concludes the proof.
\epr

\btm[MA]\label{nonmeasurable_product}
Let $G$ be a locally compact group, and let $A\subset G$ be a nonempty [locally] null set; then there is a set $S\subset G$ such that $SA$ is nonmeasurable.
\etm
\bpr
Let $H$ be an open pro-Lie subgroup of $G$. Clearly, $H$ can be chosen $\sigma$-compact (e.g., generated by any pre-compact neighborhood of identity).

Translating $A$, if necessary, we can assume that $e\in A$. Then $A_1=A\cap H$ is nonempty and [locally] null with respect to the Haar measure of $H$, and due to $\sigma$-compactness it is just null in $H$. If we find a set $S\subset H$ such that $SA_1$ is nonmeasurable in $H$, then $(SA)\cap H=SA_1$ is nonmeasurable in $G$, and so $SA$ is nonmeasurable too. We can assume therefore that $G=H$, that is: $G$ is $\sigma$-compact and pro-Lie, and $A$ is null.

%Let now $G=\liminv_{i\in I} G_i$, where every $G_i=G/K_i$ is a Lie group, $K_i$ is a compact normal subgroup of $G$, and the order on $I$ is just inclusion of $K_i$ (see remark \ref{pro-lie}). If $I$ is countable, then $G$ is Polish, and we can apply lemma \ref{polish_group}. So assume further that $I$ is uncountable.

As in the proof of Theorem \ref{rep}, either $G$ is Polish (and we can apply Lemma \ref{polish_group}), or we can find a Polish quotient $G/K$ such that $\pi(A)$ is null, where $\pi:G\to G/K$ is the quotient map. By Lemma \ref{polish_group}, there is a set $S_1\subset G/K$ such that $S_1\pi(A)$ is nonmeasurable. By Lemma \ref{quotient}, $\pi^{-1}(S_1\pi(A)) = \pi^{-1}(S_1)A$ is also nonmeasurable. Thus, we can take $S=\pi^{-1}(S_1)$, and the theorem is proved.
\epr

This theorem together with Theorem \ref{cont_hom} imply:

\btm[MA]\label{coroll_cont_hom}
Every measurable homomorphism from a locally compact group to any topological group is continuous.
\etm

In conclusion, let us review some close results. Say that a set $S$ is {\it small\,} if the union of every family of translates of $S$ of cardinality less than continuum is null. We use Martin's axiom to guarantee that every null set is small. Without MA, this depends on the set $S$. Gruenhage \cite{darji} has proved that the ternary Cantor set is small, and Darji and Keleti --- that every subset of $\R$ of packing dimension less than 1 is small. From the other side, Elekes and T\'oth \cite{elekes} and Ab\'ert~\cite{abert} proved the following: it is consistent with ZFC that in every locally compact group there is a non-small compact set of measure zero.
It is however unknown whether for a non-small set the statement of Theorem \ref{nonmeasurable_product} is false.

Finally, we say a few words on results in ZFC concerning nonmeasurable products of sets. One should better say ``sums of sets'' because there is a tradition to do everything in the commutative case. This restriction is reasonable since the principal difficulties lie already in the case of the real line. The advances most close to our topic are: for every null set $S$ on the real line such that $S+S$ has positive outer measure there is a set $A\subset S$ such that $A+A$ is nonmeasurable (Ciesielski, Fejzic and Freiling~\cite{ciesel}). Cicho\'n, Morayne, Ra\l owski, Ryll-Nardzewski, and \.Zeberski~\cite{cichon} proved that there is a subset $A$ of the Cantor set $C$ such that $A+C$ is nonmeasurable, and under additional axioms they prove the same statement for every closed null set $P$ such that $P+P$ has positive measure. There is also a series of results going back to Serpi\'nski which find null sets $A$ and $B$ such that $A+B$ is non-measurable (see, e.g., a monograph \cite{kharaz} and recent papers of Kharazishvili and Kirtadze \cite{kharaz-2004}, \cite{kharaz-2005}), where the task is to make $A$ and $B$ ``maximally negligible'' (in different senses), and $A+B$ ``maximally nonmeasurable''.

\medskip

Author thanks S.~Akbarov for drawing her attention to this problem. I thank also S.~Akbarov, J. Cicho\'n, M. Morayne, R. Ra\l owski and J. \.Zeberski for useful discussions.


\begin{thebibliography}{99}

\bibitem{abert} M. Ab\'ert. Less than continuum many translates of a compact nullset may cover any infinite profinite group. {\it J. Group Theory } {\bf 11} (2008), 545--553.
\bibitem{bourbaki} N. Bourbaki. {\it General topology, Chapters 5-10}. Springer, 1989.
\bibitem{brzdek} J. Brzd\k ek. Continuity of measurable homomorphisms. {\it Bull. Austral. Math. Soc.} {\bf 78} (2008), 171--176.
\bibitem{ryll-1979} J. Brzuchowski, J. Cicho\'n, E. Grzegorek and C. Ryll-Nardzewski, On the existence of nonmeasurable unions, {\it Bull. Polish Acad. Sci. Math. } 27 no.~6 (1979), 447--448.
\bibitem{bukovsky} L. Bukovsk\'y, Any partition into Lebesgue measure zero sets produces a non-measurable set, {\it Bull. Acad. Polon. Sci.} {\bf 27} no.~6 (1979), 431--435.
\bibitem{chris} J. P. R. Christensen, Borel structures in groups and semigroups, {\it Math. Scand.} {\bf 28} (1971), 124--128.
\bibitem{cichon} J. Cicho\'n, M. Morayne, R. Ra\l owski, C. Ryll-Nardzewski, S. \.Zeberski, On nonmeasurable unions. {\it Topology and its Applications} {\bf 154} (2007), 884--893.
\bibitem{ciesel} K. Ciesielski, H. Fejzic, C. Freiling. Measure zero sets with non-measurable sum.
{\it Real Anal. Exch.} 27 no.~2 (2001), 783--794.
\bibitem{darji} U. B. Darji, T. Keleti. Covering $\R$ with Translates of a Compact Set. {\it Proc. AMS}, {\bf 131} No.~8 (2003), 2593--2596.
%\bibitem{davis} H. F. Davis, A note on Haar measure, {\it Proc. AMS} 6, no.~2 (1955), 318--321.
\bibitem{elekes} M. Elekes, \'A. T\'oth. Covering locally compact groups by less than $2^\omega$ translates of a compact nullset. {\it Fund. Math.} {\bf 193} no.~3 (2007), 243--257.
\bibitem{engelking} R. Engelking, {\it General topology}. Taylor \& Francis, 1977.
\bibitem{fremlin-martin} D. H. Fremlin, {\it Consequences of Martin's axiom}. Cambridge tracts in mathematics, no. 84. Cambridge University Press, 1984.
\bibitem{fremlin} D. H. Fremlin. {\it Set-theoretic measure theory: Vol. V}. Torres Fremlin, 2008.
\bibitem{gaal} S. A. Gaal, {\it Linear analysis and representation theory}. Springer, 1973.
\bibitem{halmos} P. R. Halmos. {\it Measure theory}. Springer, 1974.
\bibitem{HR} E. Hewitt, K. A. Ross, {\it Abstract harmonic analysis I, II.} Springer, 3rd printing, 1997.
\bibitem{hofmor-compact} K. H. Hofmann, S. A. Morris, {\it The Structure of Compact Groups}. De Gruyter, 2006.
\bibitem{hofmor} K. H. Hofmann, S. A. Morris. {\it The Lie theory of connected pro-Lie groups}. EMS, 2007.
\bibitem{kharaz} A. Kharazishvili. {\it Nonmeasurable sets and functions}. North Holland, 2004.
\bibitem{kharaz-2004} A. Kharazishvili, A. Kirtadze. On algebraic sums of measure zero sets in uncountable commutative groups. {\it Proc. A. Razmadze Math. Inst. } {\bf 135} (2004), 97--103.
\bibitem{kharaz-2005} A. Kharazishvili. The algebraic sum of two absoultely negligible sets can be an absolutely nonmeasurable set. {\it Georgian Math. J.} {\bf 12} no.~3 (2005), 455--460.
\bibitem{kleppner} A. Kleppner. Measurable Homomorphisms of Locally Compact Groups, {\it Proc. AMS} {\bf 106}, no.~2 (1989), 391--395.
%\bibitem{lee} D.~H. Lee, Supplements for the identity component in locally compact groups. {\it Math. Z.} {\bf 104} (1968), 28--49.
\bibitem{pettis} B. J. Pettis, On continuity and openness of homomorphisms in topological groups, {\it Ann. of Math.} (2), {\bf 52} (1950), 293--308.
\bibitem{reiter} H. Reiter, J. D. Stegeman. {\it Classical harmonic analysis and locally compact groups,}
Oxford: Clarendon Press, 2000.
\bibitem{rosendal} C. Rosendal. Automatic continuity of group homomorphisms. {\it Bull. Symb. Logic} {\bf 15} no.~2 (2009), 184--214.
\bibitem{shtern} A. I. Shtern. Van der Waerden continuity theorem for semisimple Lie groups. {\it Russ. J. Math. Phys.} {\bf 13} No.~2 (2006), 210--223.
\bibitem{waerden} B. L. van der Waerden, Stetigkeitss\"atze f\"ur halbeinfache Liesche Gruppen, {\it Math. Z.} {\bf 36} no.~1 (1933), 780--786.
\end{thebibliography}
\end{document}